\documentclass[amssymb,11pt, draft]{amsart}
\usepackage[latin1]{inputenc}
\usepackage{amsmath}
\usepackage{amsfonts}

\title{Counting Free Abelian Actions}
\author{Tad White}
\address{Department of Defense}
\email{tadwhite1@gmail.com}

\newcommand{\Z}{\mathbf{Z}}
\newcommand{\Hom}{\mathrm{Hom}}
\newcommand{\Stab}{\mathrm{Stab}}

\newtheorem{theorem}{Theorem}
\newtheorem{lemma}[theorem]{Lemma}

\begin{document}

\begin{abstract}
We consider the problem of counting commuting $r$-tuples of elements of the symmetric group $S_n$, i.e.~computing $|\Hom(\Z^r,S_n)|$.  The cases $r=1,2$ are well-known; a product formula for the case $r=3$ was conjectured by Adams-Watters and later proved by Britnell.  In this note we solve the problem for arbitrary $r$.
\end{abstract}
\maketitle

\section{Introduction}

In \cite{Britnell12}, Britnell proves the following product formula for the exponential generating
function of the number $T(n)$ of ordered triples of commuting elements of the symmetric group $S_n$:
\begin{equation} \label{eq:triples}
   \sum\limits_{n=0}^{\infty} \dfrac{T(n)}{n!}\, u^n = \prod\limits_{j=1}^{\infty} (1-u^j)^{-\sigma(j)}
\end{equation}
where $\sigma(j)$ denotes the sum of the divisors of $j$.  The right-hand side of (\ref{eq:triples}) is
the so-called Euler transform of the sequence $\{\sigma(j)\}$.  The formula (\ref{eq:triples}) had been conjectured
by Adams-Watters, based on a comparison of sequences A079860 and A061256 in Sloane's Online Encyclopedia
of Integer Sequences \cite{OEIS}.

In this note we derive an analogous product formula which counts
commuting $r$-tuples in $S_n$.  As Britnell points out, the appearance of the
number-theoretic function $\sigma(j)$ in (\ref{eq:triples}) is surprising.  The present proof not only
yields a unified approach to the known results, but also illuminates the number-theoretic connections.
The proof rests on a standard combinatorial technique for enumerating
structures in terms of ``connected'' structures.

\section{The result}
For any integer $r\ge 0$, let $T_r(n)$ denote the number of $r$-tuples
$\{g_1,\ldots,g_r\}$ in $S_n$
such that $g_i g_j = g_j g_i$ for all $i$, $j$.
We begin by observing that such tuples correspond precisely to homomorphisms from the group
$\Z^r$ to $S_n$: we map the $i$-th standard basis vector to $g_i$.  That is, $T_r(n) = |\Hom(\Z^r,S_n)|$.

Equivalently, $T_r(n)$ counts the number of actions of $\Z^r$ on the $n$-set $[n]=\{1,\ldots,n\}$.
We let $T_r(n,k)$ denote the number of such actions with $k$ orbits.  For any integer $r\ge 1$, we denote
by $\lambda_r(n)$ the number of index-$n$ subgroups of $\Z^r$; for $r=0$, it will be convenient to define
$\lambda_0(n)$ to be 1 if $n=1$, and 0 otherwise.

\begin{theorem}
\label{thm:maintheorem}
For any integer $r\ge 1$, we have
\begin{equation} \label{eq:mainorbitformula}
   \sum\limits_{n,k=0}^{\infty} T_r(n,k)\frac{u^n}{n!} y^k = \prod\limits_{j=1}^{\infty} (1-u^j)^{-y\lambda_{r-1}(j)}.
\end{equation}
In particular, we can take $y=1$ to ignore the orbit counts:
\begin{equation} \label{eq:mainformula}
   \sum\limits_{n=0}^{\infty} T_r(n)\frac{u^n}{n!} = \prod\limits_{j=1}^{\infty} (1-u^j)^{-\lambda_{r-1}(j)}.
\end{equation}
\end{theorem}

It turns out that $\lambda_r(n)$ is quite tractable:

\begin{lemma}
[Counting subgroups]
\label{lemma:countingsubgroups}
For any $r\ge 1$, 
\[\lambda_r(n)= \sum\limits_{d_1d_2\cdots d_r=n}d_2 d_3^2 \cdots d_r^{r-1}.\]
Equivalently, $\{\lambda_r(n)\}$ is the Dirichlet convolution of the sequences $\{n^k\}$ for $0\le k < r$
\cite{ApostolNumThy}.
\end{lemma}

\proof Every subgroup of $\Z^r$ has a unique basis in Hermite normal form \cite{CohenCompAlgNumThy}.
If the subgroup has index $n$ in $\Z^r$, this normal form is an upper-triangular $r\times r$ matrix
in which the product of the diagonal elements $d_1,\ldots,d_r$ is $n$, and the elements in the column
above $d_i$ lie between 0 and $d_i-1$ inclusive.
Every such matrix corresponds to a unique subgroup, so the lemma follows by counting the possible
Hermite normal forms.\qed

We recover some known results when $r$ is small:

\begin{enumerate}
\item[$r=1$:]  The right-hand side of (\ref{eq:mainformula}) is just
$1/(1-u) = \sum_{n\ge 0}u^n$, in accordance with the trivial fact that $T_1(n)=|S_n|=n!$. 

\item[$r=2$:]  $\lambda_1(j)=1$, as $\Z$ has a unique subgroup of each index $j$.  So the right-hand side of (\ref{eq:mainformula})
becomes \[\prod_{j\ge 1} (1-u^j)^{-1} = \sum_{n\ge 0}p(n) u^n,\] where $p(n)$ is the number of integer
partitions of $n$.  Hence $T_2(n)= p(n)\,n!$.  (Britnell points out that this result appears in
Erd\H{o}s and Tur\'{a}n \cite{ErdosTuran}; it can also be derived easily by applying Burnside's lemma
to the action of $G$ on itself by conjugation.)

\item[$r=3$:]  This is the case considered by Britnell.  By Lemma \ref{lemma:countingsubgroups},
$\lambda_2(j)=\sum_{d|j} d = \sigma(j)$, so (\ref{eq:triples}) follows.
\end{enumerate}

In case $r=4$, we have
\begin{equation}
\label{eq:lambda4}
  \lambda_3(j)=\sum_{d e|j} d e^2 = \sum_{d|n} d^2 \sigma(n/d),
\end{equation}
which is A001001 in \cite{OEIS}.  $T_4(n)/n!$ is the Euler transform of this sequence, namely
$\{1, 1, 8, 21, 84, 206,\ldots\}$.

\section{The proof}
\label{sec:Proof}

One frequently encounters objects which can be decomposed into a number of independent ``connected'' objects.
The so-called exponential formula \cite{gfology} is a general combinatorial principle asserting that the
exponential generating function (egf) counting the labeled objects is the exponential of the egf
counting the labeled connected objects.

In the present context,
we wish to count actions of a group $G$ on a set $X$.  The appropriate notion of ``connected'' in this
context is ``transitive,'' since $X$ is the union of its (connected) $G$-orbits.
(Put another way, we can associate
to each action a graph with vertex set $X$, joining two vertices by an edge if some group element sends
one to the other.  The connected components of this graph are just the orbits of the action; the graph
is connected iff the action is transitive.)

Let $d_n$ denote the number of transitive actions of a group $G$ on $[n]$, and let
\[\mathcal{D}(u) = \sum\limits_{n=0}^{\infty} d_n \frac{u^n}{n!}\]
be the egf of $\{d_n\}$.  Let $h(n,k)$ denote the number of
$G$-actions on $[n]$ containing $k$ orbits.  The exponential formula asserts in this case that:
\begin{equation}
\label{eq:ExpFormula}
\sum\limits_{n,k=0}^{\infty} h(n,k) \frac{u^n}{n!} y^k = \mathrm{exp}\left(y \mathcal{D}(u)\right)
\end{equation}
A direct derivation of this result can be found in Lubotzky \cite[Prop.~1.10]{LubotzkyCountingSubgroups}.

Our task is thus to compute $\mathcal{D}(u)$ when $G=\Z^r$, i.e.~to enumerate the transitive actions of $\Z^r$ on $[n]$.
We begin by observing that transitive actions have a very simple form (see also \cite[Prop.~1.1]{LubotzkyCountingSubgroups}):

\begin{lemma}[Transitive actions are coset actions]
\label{lemma:cosetactions}
Let $G$ act transitively on a set $X$.  Then there exists a unique subgroup $K\subset G$ such
that the given action is equivalent to the action of $G$ on the cosets $K\backslash G$;
that is, there is a bijection $f:X\to K\backslash G$ such that
\begin{equation}
\label{eq:bijection} f(x) g = f(x g)\textrm{\quad for all $x\in X$, $g\in G$.}
\end{equation}
Furthermore, there are precisely $|X|$ such bijections.
\end{lemma}

\proof
Choose a basepoint $\ast\in X$, and let $K=\Stab(\ast)$.  By transitivity, for any $x\in X$,
we can select $g_x\in G$ such that $\ast g_x = x$.  Note that if $g$ is another such choice,
we have $\ast g_x g^{-1}=\ast$, so $K g_x = K g.$   Thus the map $f: X\to K\backslash G$ given by
$x\mapsto K g_x$ is well-defined (independent of $g_x$), and (\ref{eq:bijection}) clearly holds.

The $|X|$ bijections arise from the $|X|$ possible choices of $\ast$: once $f(\ast)$ has been chosen,
(\ref{eq:bijection}) determines $f$ completely by transitivity. 
\qed

It follows that we can count transitive actions of $G=\Z^r$ on $[n]$ by counting the possible point
stabilizers $K$.  $K$ must be an index-$n$ subgroup of $\Z^r$.
Given such a $K$, and the standard action of $\Z^r$ on $\Z^r/K$, we can use any bijection
$f: [n] \to \Z^r/K$ to transfer the action to one on $[n]$.  By Lemma \ref{lemma:cosetactions}, each
action on $[n]$ arises in exactly $n$ ways.  Hence each possible $K$ corresponds to
$n!/n = (n-1)!$ actions.  As there are $\lambda_r(n)$ choices for $K$, we see that
the number of transitive actions of $\Z^r$ on $[n]$ is $(n-1)!\,\lambda_r(n)$.

We can now compute, for fixed $r>0$, the egf $\mathcal{D}(u)$ for transitive actions of $\Z^r$ on $[n]$:

\begin{eqnarray}
\mathcal{D}(u)&=&
\sum\limits_{n=0}^{\infty} (n-1)!\,\lambda_r(n)\frac{u^n}{n!}\nonumber\\
 &=&\sum\limits_{n=0}^{\infty} \lambda_r(n)\frac{u^n}{n}\nonumber\\
 &=& \sum\limits_{n=0}^{\infty} \frac{u^n}{n} \sum\limits_{d_1d_2\cdots d_r=n} d_2\, d_3^2 \cdots d_r^{r-1}
     \quad\textrm{by Lemma \ref{lemma:countingsubgroups}}\nonumber\\
 &=& \sum\limits_{n=0}^{\infty} u^n \sum\limits_{\ d_1d_2\cdots d_r=n} \frac{1}{d_1}\, d_3 d_4^2 \cdots d_r^{r-2}\nonumber\\
 &=& \sum\limits_{d_1, d_2, \ldots, d_r} d_3 d_4^2 \cdots d_r^{r-2} \frac{u^{d_1\cdots d_r}}{d_1}\nonumber\\
 &=& \sum_{j\ge 1}-\lambda_{r-1}(j)\log(1-u^{j}) \quad\textrm{after taking $j=d_2\cdots d_r$.}\nonumber
\end{eqnarray}

Theorem \ref{thm:maintheorem} now follows immediately from the exponential formula.\qed
\relax\par
\medskip
We remark that, while we have considered these power series formally, everything converges within the
unit disk. From Lemma \ref{lemma:countingsubgroups}, since each $d_i$ is at most $n$, we have
$\lambda_r(n)<n^{r(r-1)/2}$,
which is polynomial in $n$.  Thus $\sum \lambda_r(n) u^n/n$ converges in the unit disk, and therefore
so do the series on the left-hand sides of (\ref{eq:mainorbitformula}) and (\ref{eq:mainformula}).


\bibliographystyle{plain}
\bibliography{myRefs}

\end{document}